\documentclass{amsart}
\title[On a theta product of Jacobi...]
{On a theta product of Jacobi and its applications to $q$-gamma products}
\usepackage{amssymb,amsmath,amsthm,epsfig,graphics,latexsym}
 \theoremstyle{definition}

  \theoremstyle{plain}

  \newtheorem{theorem}    {Theorem}

  \theoremstyle{remark}
  \newtheorem{example}{{\bf Example}}
  \newtheorem{rmk}{{\bf Remark}}
  \newcommand{\te}{\theta}
  \newcommand{\fr}{\frac}
  \newcommand{\lt}{\left(}
  \newcommand{\rt}{\right)}
\begin{document}
  \author[M. El Bachraoui and J. S\'{a}ndor]{Mohamed El Bachraoui and J\'{o}zsef S\'{a}ndor}
  \address{Dept. Math. Sci,
 United Arab Emirates University, PO Box 15551, Al-Ain, UAE}
 \email{melbachraoui@uaeu.ac.ae}
 \address{Babes-Bolyai University, Department of Mathematics and Computer Science, 400084 Cluj-Napoca, Romania}
 \email{jsandor@math.ubbcluj.ro}
 \begin{abstract}
 We give a new proof for a product formula of Jacobi which turns out to be equivalent to a
$q$-trigonometric product which was stated without proof by
Gosper. We apply this formula to derive a $q$-analogue for the Gauss multiplication formula for
the gamma function. Furthermore, we give explicit formulas for short products of $q$-gamma functions.
 \end{abstract}
  \date{\textit{\today}}
  \maketitle
\section{Introduction}
Throughout we let $\tau$ be a complex number in the upper half plane, let $q=e^{\pi i\tau}$, and let $\tau'=-\frac{1}{\tau}$.
Note that the assumption $\mathrm{Im}(\tau)>0$ guarantees
that $|q|<1$.
The $q$-shifted factorials of a complex number $a$ are defined by
\[
(a;q)_0= 1,\quad (a;q)_n = \prod_{i=0}^{n-1}(1-a q^i),\quad
(a;q)_{\infty} = \lim_{n\to\infty}(a;q)_n.
\]
It is easily verified that
\begin{equation}\label{q-basics}
(q^2;q^2)_{\infty} = (q;q)_{\infty} (-q;q)_{\infty},\quad
(q;q)_{\infty} = (q^2;q^2)_{\infty} (q;q^2)_{\infty},
\end{equation}
\[
(-q;q)_{\infty}  = \fr{1}{(q;q^2)_{\infty}}, \quad\text{and\quad }
(-q;q)_{\infty} = (-q^2;q^2)_{\infty} (-q;q^2)_{\infty}.
\]
Ramanujan theta functions $\psi(q)$ and $f(q)$ are given by
\[
\psi(q) = \fr{(q^2;q^2)_{\infty}}{(q;q^2)_{\infty}}\quad\text{and\quad}
f(-q) = (q;q)_{\infty}.
\]
See Berndt~\cite{Berndt-1} for some properties of Ramanujan theta functions.
For convenience we write
\[
(a_1,\ldots,a_k;q)_n = (a_1;q)_n\cdots (a_k;q)_n,\quad
(a_1,\ldots,a_k;q)_{\infty} = (a_1;q)_{\infty} \cdots (a_k;q)_{\infty}.
\]
Jacobi first theta function is defined as follows:
\[
\theta_1(z,q)  =\theta_1(z \mid \tau) = 2\sum_{n=0}^{\infty}(-1)^n q^{(2n+1)^2/4}\sin(2n+1)z.
\]
 The numbers $\tau$ and $q$ respectively are referred to as the \emph{parameter} and \emph{nome}
of the theta functions.
An important property of Jacobi theta functions is their infinite product
representations which for the function $\theta_1(z\mid\tau)$ is given by
\begin{equation}\label{theta-product}
\theta_1(z \mid\tau) =
i q^{\frac{1}{4}}e^{-iz} (q^2 e^{-2iz},e^{2iz},q^2; q^2)_{\infty} .
\end{equation}
Jacobi~\cite{Jacobi} proved that
\begin{equation}\label{MainProd}
\fr{(q^{2n};q^{2n})_{\infty}}{(q^2;q^2)_{\infty}^n}
\prod_{k=-\fr{n-1}{2}}^{\fr{n-1}{2}}\te_1 \left(z+\fr{k\pi}{n} \bigm| \tau \right) =
\te_1(nz \mid n\tau),
\end{equation}
see also Enneper~\cite[p. 249]{Enneper}.
Unlike many of Jacobi's results the formula (\ref{MainProd})
seems not to have received much attention by mathematicians. This is probably due to the lack of applications.
Furthermore, up to the authors' knowledge no new proof has been given for this product formula.
Besides, this formula turns out to be equivalent to a $q$-trigonometric identity
of Gosper (see (\ref{SineProd}) below) which he apparently was not aware of as
he stated the identity without mentioning any reference to it.
Our first goal in this note is to offer a new proof for (\ref{MainProd}).
To this end we will need the following basic properties of the function $\theta_1(z\mid\tau)$.
\begin{align}\label{theta-1}
\theta_1(k\pi\mid\tau) &= 0 \quad (k\in\mathbb{Z}), \nonumber \\
\theta_1(z+\pi \mid \tau) & = \te_1(-z\mid \tau) = -\te_1(z\mid \tau),  \\
\theta_1(z+\pi\tau \mid \tau) &= -q^{-1} e^{-2iz} \theta_1(z \mid \tau). \nonumber
\end{align}
It can be shown that the last formula can be extended as follows
 \begin{equation} \label{Transf}
 \te_1 \left( z+\pi\tau \bigm| \frac{\tau}{k} \right) = (-1)^k q^{-k}e^{-2kiz}
 \te_1 \left(z\bigm| \frac{\tau}{k} \right)\quad (k\in\mathbb{N}).
 \end{equation}
  We will need
Jacobi's imaginary transformation stating that
\begin{equation}\label{ImTransf}
\theta_1(z\mid \tau')
=-i (-i\tau)^{\frac{1}{2}}  e^{\frac{i\tau z^2}{\pi}}\theta_1(z\tau \mid \tau).
\end{equation}
Letting $\te'(z\mid\tau)$ denote the first derivative of $\theta(z\mid\tau)$ with respect to $z$, we have
\begin{equation}\label{theta-1-deriv}
\te_1'(0\mid\tau) = 2 q^{\fr{1}{4}} (q^2;q^2)^3
\end{equation}
and
\begin{equation} \label{key-derivative}
\te_1'(0\mid \tau') = 2 (-i\tau)^{\fr{3}{2}}  q^{\fr{1}{4}} (q^2;q^2)^3.
\end{equation}
For details about theta functions we refer to the book by Whittaker~and~Watson~\cite{Whittaker-Watson} and the book by Lawden~\cite{Lawden}. For recent references on Jacobi theta
 functions which are closely related to our current topic, the reader is referred to
 Liu~\cite{Liu-2005, Liu-2007} and Shen~\cite{Shen-1, Shen-2}.
As to our applications, we shall use
(\ref{MainProd}) in an equivalent form to establish new $q$-analogues for well-known products involving the gamma function
$\Gamma (z)$ as we describe now.
The $q$-gamma function is given by
\[
\Gamma_q(z) = \dfrac{(q;q)_\infty}{(q^{z};q)_\infty} (1-q)^{1-z}  \quad (|q|<1).
\]
It is immediate from the previous definition and (\ref{q-basics}) that
\begin{equation}\label{q-gamma-half}
\Gamma_q\left(\fr{1}{2}\right) = \fr{(q;q)_{\infty}}{(q^{\fr{1}{2}};q)_{\infty}} \sqrt{1-q}
= \fr{f^2(-q)}{f(-q^{\fr{1}{2}})} \sqrt{1-q} = \psi(q^{\fr{1}{2}})\sqrt{1-q}.
\end{equation}
It is well-known that $\Gamma_q (z)$ is a $q$-analogue for the function $\Gamma (z)$, see Gasper~and~Rahman~\cite{Gasper-Rahman}.

\noindent
The Gaussian multiplication formula for the gamma function states that
\begin{equation}\label{Gauss-Gamma-special}
\Gamma\left(\fr{1}{n}\right) \Gamma\left(\fr{2}{n}\right) \cdots \Gamma\left(\fr{n-1}{n}\right) =
 \fr{(2\pi)^{\fr{n-1}{2}}}{\sqrt{n}} \quad (n=1,2,\ldots).
\end{equation}
A natural question is if, for the product
\begin{equation}\label{Sandor-problem}
\Gamma_q\left(\fr{1}{n}\right) \Gamma_q\left(\fr{2}{n}\right) \cdots \Gamma_q\left(\fr{n-1}{n}\right),
\end{equation}
one can find some closed formula of type (\ref{Gauss-Gamma-special}).
We note that~(\ref{Gauss-Gamma-special}) can be deduced by an application of the following well-known more general identity:
 \begin{equation}\label{Gauss-Gamma}
 n^{nz-\fr{1}{2}} \Gamma(z)\Gamma\left(z+\fr{1}{n}\right) \cdots \Gamma\left(z+\fr{n-1}{n}\right) =
 \Gamma(nz) (2\pi)^{\fr{n-1}{2}}\quad (n=1,2,\ldots).
 \end{equation}
 A famous $q$-analogue for (\ref{Gauss-Gamma}) due to Jackson~\cite{Jackson-1, Jackson-2},
 (see \cite[p. 22]{Gasper-Rahman}), states that
 \begin{equation}\label{Jackson-q-Gamma}
 \left(\fr{1-q^n}{1-q}\right)^{nz-1} \Gamma_{q^n} (z)\Gamma_{q^n} \left(z+\fr{1}{n}\right)
 \cdots \Gamma_{q^n} \left(z+\fr{n-1}{n}\right)
 \end{equation}
 \[
 =
 \Gamma_{q^n} (nz)\Gamma_{q^n} \left(\fr{1}{n}\right) \Gamma_{q^n} \left(\fr{2}{n}\right)
 \cdots \Gamma_{q^n} \left(\fr{n-1}{n}\right) \quad (n=1,2,\ldots).
 \]
 However, it seems not to be easy to derive a closed
 formula for the product~(\ref{Sandor-problem}) using the relation~(\ref{Jackson-q-Gamma}).
 Our secondary goal in this note is to apply Jacobi's relation~(\ref{MainProd}) to establish
 a closed formula for the product~~(\ref{Sandor-problem}). This gives a $q$-analogue for the
 formula~(\ref{Gauss-Gamma-special}) which seems not to appear in literature.

 Furthermore, S\'{a}ndor~and~T\'{o}th~\cite{Sandor-Toth} found
 \begin{equation}\label{short-gam-prod}
 P(n):= \prod_{\substack{k=1 \\ (k,n)=1}}^n \Gamma\left(\fr{k}{n}\right)
 = \fr{(2\pi)^{\fr{\varphi(n)}{2}}}{ \prod_{d\mid n} d^{\fr{1}{2}\mu\left(\fr{n}{d}\right)} }
 = \fr{(2\pi)^{\fr{\varphi(n)}{2}}}{ e^{\fr{\Lambda(n)}{2}} },
 \end{equation}
 where $\varphi(n)$ is the Euler totient function, $\mu(n)$ is the M\"{o}bius mu function, and $\Lambda(n)$ is the Von~Mangoldt function.
 We accordingly let
 \[
 P_q(n) = \prod_{\substack{k=1 \\ (k,n)=1}}^n \Gamma_q\Big(\fr{k}{n}\Big).
 \]
 Our third purpose in this note is to evaluate the last product and therefore give a $q$-version of the short product (\ref{short-gam-prod}).
 To have our formula look like (\ref{short-gam-prod}) we introduce the $q$-Von~Mangoldt function as follows:
 \[
 \Lambda_q(n)=
 \log \fr{ 2^{\varphi(n)}\prod_{d\mid n} \big(f(-q^{\fr{1}{d}})\big)^{2\mu\left(\fr{n}{d}\right)} }{(q^{\fr{1}{2}};q)_{\infty}^{2\varphi(n)} }.
 \]
It turns out that our formula for $P_q(n)$ when $n$ is a power of $2$ can be expressed in terms of
Ramanujan function $\psi$. This and some work by Berndt~\cite{Berndt-2},
Yi~\emph{et al.}~\cite{Yi-Lee-Paek}, and Baruah and Saikia~\cite{Baruah-Saikia} enable us to deduce
explicit identities for a variety of short products of $q$-gamma functions.
 For references on short products of the gamma function we refer to \cite{BenAri-et-al, Chamberland-Straub, Martin, Nijenhuis, Nimbran}.

\noindent
To derive our results on products of the $q$-gamma function we shall use the link of this function with
the $q$-trigonometry of Gosper.
Gosper~\cite{Gosper} introduced $q$-analogues of $\sin z$ and $\cos z$ as follows
\begin{equation}\label{sine-cosine-q-gamma}
\begin{split}
\sin_q \pi z
&=
q^{\fr{1}{4}} \Gamma_{q^2}^2\left(\fr{1}{2}\right) \frac{q^{z(z-1)}}{\Gamma_{q^2}(z) \Gamma_{q^2}(1-z)} \\
\cos_q \pi z
&=
\Gamma_{q^2}^2\left(\fr{1}{2}\right) \fr{q^{z^2}}{\Gamma_{q^2}\left(\fr{1}{2}-z \right) \Gamma_{q^2}\left(\fr{1}{2}+z\right)}.
\end{split}
\end{equation}
It can be shown that
$\lim_{q\to 1}\sin_q z = \sin z$ and $\lim_{q\to 1}\cos_q z = \cos z$.
Gosper proved that the functions $\sin_q$ and $\cos_q z$ are related to the function
$\theta_1(z\mid \tau')$ as follows:
\begin{equation}\label{sine-cosine-theta}
\sin_q (z) = \frac{\theta_1(z\mid \tau')}{\theta_1\left( \frac{\pi}{2}\bigm| \tau' \right)} \qquad \text{and \quad}
\cos_q (z) = \frac{\theta_1\left( z+\frac{\pi}{2} \bigm| \tau' \right)}
{\theta_1 \left( \frac{\pi}{2} \bigm| \tau' \right)} \quad \quad (\tau' = \frac{-1}{\tau})
\end{equation}
from which it immediately follows that $\sin_q (z-\pi/2)=\cos_q z$, $\sin_q \pi = 0$,
and $\sin_q \fr{\pi}{2} = 1$.
Gosper, on the one hand, stated many identities involving $\sin_q z$ and $\cos_q z$ which easily
follow just
from the definition and the basic properties of the function $\theta_1(z\mid\tau)$. For instance, he derived that
\begin{equation}\label{q-sin-derive}
\sin_q' 0 = \fr{-2 \ln q}{\pi} q^{\fr{1}{4}} \fr{(q^2;q^2)_{\infty}^2}{(q;q^2)_{\infty}^2} =
-\fr{2\ln g}{\pi}q^{\fr{1}{4}}\psi^2(q).
\end{equation}
On the other hand, Gosper~\cite{Gosper} using the computer facility \emph{MACSYMA} stated without proof a variety of identities involving $\sin_q z$ and $\cos_q z$ and he
asked the natural question whether his formulas hold true.
For recent work on Gosper's conjectures we refer to \cite{Touk-Houchan-Bachraoui, Bachraoui-1, Bachraoui-2, Bachraoui-3, Mezo-1}.
Among the formulas which Gosper~\cite[p. 92]{Gosper} stated without proof we have
\begin{equation}\label{SineProd}
\prod_{k=0}^{n-1}\sin_{q^n}\pi \left(z+\fr{k}{n} \right) =
q^{\frac{(n-1)(n+1)}{12}} \frac{(q;q^2)_{\infty}^2}{(q^n;q^{2n})_{\infty}^{2n}} \sin_q n\pi z.
\end{equation}
 However, by using the relation (\ref{sine-cosine-theta}) and some basic manipulations, one can show
 that (\ref{SineProd}) is actually equivalent to Jacobi's multiplication~(\ref{MainProd}).
 \section{Main results and some examples}
We start with results on the $q$-gamma function.
 \begin{theorem}\label{Gauss-q-Gamma}
 For any positive integer $n$, there holds
 \begin{equation}\label{q-anlog-Gauss}
\prod_{k=1}^{n-1}\Gamma_q\left(\fr{k}{n}\right) =
\Big(\Gamma_q\big(\fr{1}{2}\big)\Big)^{n-1} \fr{f^{n-1}(-q^{\fr{1}{2}})}{f^{n-2}(-q) f(-q^{\fr{1}{n}})}.
\end{equation}
Moreover, identity (\ref{q-anlog-Gauss}) is the $q$-analogue for identity~(\ref{Gauss-Gamma-special}).
 \end{theorem}
\begin{rmk}\label{rmk:general-prod}
We note that Mahmoud and Agarwal in~\cite[Theorem 7]{Mahmoud-Agarwal} proved that for
$x>0$ and $0<q<1$
\begin{equation}\label{M-A}
\prod_{k=0}^{n-1}\Gamma_{q^n}\Big(\fr{x+k}{n}\Big)  = \fr{(q^n;q^n)_{\infty}^n(1-q)^{\fr{n-1}{2}}}{(q;q)_{\infty}}
\Big(\fr{1-q^n}{1-q}\Big)^{1-x} \Gamma_q (x).
\end{equation}
However, their formula is incorrect as on its right hand-side there should be a factor
$(1-q^n)^{\fr{n-1}{2}}$ instead of the factor $(1-q)^{\fr{n-1}{2}}$. 
Moreover, their proof is complicated. We improve the authors' formula by means
of our arguments as follows.
Combining (\ref{Jackson-q-Gamma}) and Theorem~\ref{Gauss-q-Gamma} yields
\begin{equation}\label{rmk-help-1}
\prod_{k=0}^{n-1}\Gamma_{q^n}\Big(z+\fr{k}{n}\Big)
= \fr{(q^n;q^n)_{\infty}^n(1-q^n)^{\fr{n-1}{2}}}{(q;q)_{\infty}}
\Big(\fr{1-q^n}{1-q}\Big)^{1-nz} \Gamma_q (nz),
\end{equation}
which by letting $z=\fr{x}{n}$ improves (\ref{M-A}). Furthermore, putting $q^{\fr{1}{n}}$ in place of $q$ in
(\ref{rmk-help-1}) gives
\begin{equation}\label{rmk-help-2}
\prod_{k=0}^{n-1}\Gamma_{q}\Big(z+\fr{k}{n}\Big)
= \fr{(q;q)_{\infty}^n(1-q)^{\fr{n-1}{2}}}{(q^{\fr{1}{n}};q^{\fr{1}{n}})_{\infty}}
\Big(\fr{1-q}{1-q^{\fr{1}{n}}}\Big)^{1-nz} \Gamma_{q^{\fr{1}{n}}} (nz),
\end{equation}
which clearly extends Theorem~\ref{Gauss-q-Gamma} by letting $z=\fr{1}{n}$.
\end{rmk}
\begin{theorem}\label{q-short-prod}
For any positive integer $n$, there holds
\[
P_q(n) = \fr{ \left(\Gamma_q\Big( \fr{1}{2} \Big)\right)^{\varphi(n)} (q^{\fr{1}{2}};q)^{\varphi(n)}}
{ \prod_{d\mid n}\big( f(- q^{\fr{1}{d}}\big) \big)^{\mu\left(\fr{n}{d}\right)}} 
= \fr{\left( 2\Gamma_q\Big(\fr{1}{2}\Big)\right)^{\fr{\varphi(n)}{2}}}{e^{\fr{\Lambda_q(n)}{2}}}.
\]
\end{theorem}
Note that $\lim_{q\to 1} P_q(n)  = P(n)$ and so
$\lim_{q\to 1} \Lambda_q(n) = \Lambda (n)$.
%
%
We shall now provide examples of explicit values for some products of $q$-gamma functions.
To this end, we will use
some results of Berndt~\cite{Berndt-2}, Yi~{\em et al.}~\cite{Yi-Lee-Paek},
and Baruah and Saikia~\cite{Baruah-Saikia} on explicit identities of Ramanujan's $\psi$ function.
We start by a more general result.
\begin{theorem}\label{2-powers}
For any positive integer $m>1$, there holds
\[
P_q(2^m) = \prod_{k=1}^{2^{m-1}}\Gamma_q\Big(\fr{2k-1}{2^m}\Big)
= (1-q)^{2^{m-2}} \psi(q^{\fr{1}{2^{m}}}) \prod_{k=1}^{m-1} \psi^{2^{m-1-k}}(q^{\fr{1}{2^{k}}}).
\]
\end{theorem}
We now list the explicit values of the $\psi$ function which are needed for our goal. Throughout this section let
$a=\fr{\pi^{1/4}}{\Gamma\big(\fr{3}{4}\big)}$.
The following are due to Berndt~\cite[p. 325]{Berndt-2}
\[
\begin{split}
\psi(e^{-\pi}) &= a 2^{-5/8} e^{\pi/8}, \\
\psi(e^{-2\pi}) &= a 2^{-5/4} e^{\pi/4}, \\
\psi(e^{-\pi/2}) &= a 2^{-7/16}(\sqrt{2}+1)^{1/4} e^{\pi/16},
\end{split}
\]
the following are found by Yi~\emph{et al.}~\cite{Yi-Lee-Paek}
\[
\psi(-e^{-\pi}) = a 2^{-3/4} e^{\pi/8} \quad\text{and\quad }
\psi(-e^{-2\pi})= a 2^{-15/16} e^{\pi/4},
\]
and the following is given by Baruah and Saikia~\cite{Baruah-Saikia}
\[
\psi(-e^{-\pi/2}) = a 2^{-7/16} e^{\pi/16} (\sqrt{2}-1)^{1/4}.
\]
We are now ready to produce some concrete examples.
\begin{example}
Letting in Theorem~\ref{2-powers}, $m=2$, we obtain
\[
\begin{split}
\Gamma_{e^{-2\pi}} \Big(\fr{1}{4}\Big) \Gamma_{e^{-2\pi}}\Big(\fr{3}{4}\Big)
&=
(1-e^{-2\pi}) \psi(e^{-\pi}) \psi(e^{-\pi/2}) \\
&= (1-e^{-2\pi}) a^2 2^{-17/16} e^{3\pi/16} (\sqrt{2}+1)^{1/4}, \\
\Gamma_{e^{-4\pi}} \Big(\fr{1}{4}\Big) \Gamma_{e^{-4\pi}}\Big(\fr{3}{4}\Big)
&=
(1-e^{-4\pi}) \psi(e^{-2\pi}) \psi(e^{-\pi}) \\
&= (1-e^{-4\pi}) a^2 2^{-15/8} e^{3\pi/8}, \\
\Gamma_{-e^{-2\pi}} \Big(\fr{1}{4}\Big) \Gamma_{-e^{-2\pi}}\Big(\fr{3}{4}\Big)
&=
(1+e^{-2\pi}) \psi(-e^{-\pi}) \psi(-e^{-\pi/2}) \\
&=
(1+e^{-2\pi}) a^2 2^{-19/16} e^{3\pi/16}, \\
\Gamma_{-e^{4\pi}} \Big(\fr{1}{4}\Big) \Gamma_{-e^{-2\pi}}\Big(\fr{3}{4}\Big)
&=
(1+e^{-4\pi}) \psi(-e^{-2\pi}) \psi(-e^{-\pi}) \\
&=
(1+e^{-4\pi}) a^2 2^{-27/16} e^{3\pi/8}.
\end{split}
\]
Note that the first two identities in the previous list were first obtained by Mez\H{o}~\cite{Mezo-2}.
\end{example}
\begin{example}
Let $m=3$ in Theorem~\ref{2-powers}. Then
\[
\begin{split}
\Gamma_{e^{-4\pi}} \Big(\fr{1}{8}\Big) \Gamma_{e^{-4\pi}}\Big(\fr{3}{8}\Big) \Gamma_{e^{-4\pi}} \Big(\fr{5}{8}\Big) \Gamma_{e^{-4\pi}}\Big(\fr{7}{8}\Big)
&= (1-e^{-4\pi})^2 \psi^2(e^{-2\pi}) \psi(e^{-\pi}) \psi(e^{-\pi/2}) \\
&= (1-e^{-4\pi})^2 a^4 2^{-57/16} e^{11\pi/16} (\sqrt{2}+1)^{1/4},
\end{split}
\]
and similarly
\[
\Gamma_{-e^{-4\pi}} \Big(\fr{1}{8}\Big) \Gamma_{-e^{-4\pi}}\Big(\fr{3}{8}\Big) \Gamma_{-e^{-4\pi}} \Big(\fr{5}{8}\Big) \Gamma_{-e^{-4\pi}}\Big(\fr{7}{8}\Big)
=
(1+e^{-4\pi})^2 a^4 2^{-49/16} e^{11\pi/16}.
\]
\end{example}
\section{A new proof for  Jacobi's identity (\ref{MainProd})}
\noindent
We shall prove (\ref{SineProd}) which as noticed before is an equivalent form of (\ref{MainProd}).
We will employ the following result.
\begin{theorem}\cite{Bachraoui-3} \label{master}
Let $n$ be a positive integer and let $f(u)$ be an entire function such that
\[
f(u+\pi) = - f(u)\quad \text{and\ } f\lt u+\fr{\pi\tau}{n} \rt = (-1)^n q^{\fr{-1}{n}} e^{-2iu} f(u).
\]
Then for all complex numbers $x_1, x_2, \ldots, x_{n+1}$  we have:
\begin{equation*}
\sum_{j=1}^{n+1} \fr{ \te_1 \big((n-1)x_j - x_1-x_2-\cdots - x_{j-1} - x_{j+1}-x_{j+2}- \cdots - x_{n+1} \mid \tau \big) f(x_j)}
{ {\displaystyle\prod_{\substack{k=1\\ k\not= j}}^n} \te_1 \lt x_j-x_k \bigm| \fr{\tau}{n} \rt} = 0.
\end{equation*}
\end{theorem}
Note that by virtue of (\ref{sine-cosine-theta}), we can check that the desired formula (\ref{SineProd}) is equivalent to
\[
\theta_1 \big(z\mid \frac{\tau'}{n} \big) \theta_1 \left(z+\frac{\pi}{n} \bigm| \frac{\tau'}{n} \right)
\theta_1 \left(z+\frac{2 \pi}{n} \bigm| \frac{\tau'}{n} \right) \cdots
\theta_1 \left(z+\frac{(n-1)\pi}{n} \bigm| \frac{\tau'}{n} \right)
\]
\begin{equation}\label{Equiv-SineProd}
= q^{\frac{(n-1)(n+1)}{12}} \frac{(q;q^2)_{\infty}^2}{(q^n;q^{2n})_{\infty}^{2n}}
\frac{\theta_1^n \big( \frac{\pi}{2}\mid \frac{\tau'}{n} \big)}{\theta_1 \big(\frac{\pi}{2}\mid \tau')}
 \theta_1(nz \mid \tau').
\end{equation}
Next observe that the sum in Theorem~\ref{master} is equivalent to
\begin{align}\label{help-prod-0}
\te_1\big( (n-1)x_1-x_2-\ldots-x_{n+1}\mid\tau'\big) f(x_1) & \prod_{\substack{1\not=j\leq n\\ j<k\leq n+1}}\te_1(x_j-x_k\bigm| \fr{\tau'}{n})
\nonumber \\
-
\te_1\big( (n-1)x_2-x_3-\ldots-x_{n+1}-x_1 \mid\tau'  \big) f(x_2) & \prod_{\substack{2\not=j\leq n\\ j<k\leq n+1}}\te_1(x_j-x_k\bigm| \fr{\tau'}{n})
\nonumber \\
+
\te_1\big( (n-1)x_3-\ldots-x_{n+1}-x_1-x_2\mid\tau'  \big) f(x_3)& \prod_{\substack{3\not=j\leq n\\ j<k\leq n+1}} \te_1(x_j-x_k\bigm| \fr{\tau'}{n})
 \\
+ \ldots + (-1)^n \te_1\big( (n-1)x_{n+1}-x_1-\ldots-x_n\mid\tau'  \big) f(x_{n+1}) & \prod_{\substack{n+1\not=j\leq n\\ j<k\leq n+1}}\te_1(x_j-x_k\bigm| \fr{\tau'}{n})
\nonumber \\
= 0. \nonumber
\end{align}
Let $(n-1) x_{n+1} = x_1+x_2+\ldots + x_n$. Then the last term in (\ref{help-prod-0}) vanishes and for all $j=1,\ldots, n$
\[
(n-1)x_j - x_{j+1}-x_{j+2}-\ldots - x_{n+1}-x_1-\ldots - x_{j-1} = n x_j - (x_1+\ldots + x_n) - x_{n+1} = n (x_j-x_{n+1}).
\]
Then (\ref{help-prod-0}) becomes
\begin{align}\label{help-prod-1}
\te_1\big( n(x_1-x_{n+1})\mid\tau' \big) f(x_1) & \prod_{\substack{1\not=j\leq n\\ j<k\leq n+1}}\te_1(x_j-x_k\bigm| \fr{\tau'}{n}) \nonumber
\\
-
\te_1 \big( n(x_2-x_{n+1})\mid\tau' \big) f(x_2)& \prod_{\substack{2\not=j\leq n\\ j<k\leq n+1}}\te_1(x_j-x_k\bigm| \fr{\tau'}{n}) \nonumber
\\
+
\te_1 \big( n(x_3-x_{n+1})\mid\tau' \big) f(x_3) & \prod_{\substack{3\not=j\leq n\\ j<k\leq n+1}} \te_1(x_j-x_k\bigm| \fr{\tau'}{n})
\\
+ \ldots
+ (-1)^{n-1}\te_1 \big( n(x_n-x_{n+1})\mid\tau' \big) f(x_n) & \prod_{\substack{n\not=j\leq n\\ j<k\leq n+1}}\te_1(x_j-x_k\bigm| \fr{\tau'}{n})
\nonumber \\
= 0. \nonumber
\end{align}
Now assume, for $3\leq k\leq n+1$, that
\[
x_k-x_3 = \fr{(3-k)\pi}{n}.
\]
Then for all $3\leq j<k \leq n+1$
\[
\te_1 \big( n(x_k - x_j) \mid\tau' \big) = \te_1\big( (j-k)\pi \mid\tau' \big) = 0.
\]
Thus, the formula (\ref{help-prod-1}) after some simplification boils down to
\begin{equation*}
\begin{split}
\te_1 \big(n(x_1 -x_3)+(n-2)\pi \mid\tau' \big) & f(x_1)  \te_1\left( x_2-x_3 \bigm| \fr{\tau'}{n} \right)  \te_1 \left( x_2-x_3+ \fr{\pi}{n} \bigm| \fr{\tau'}{n} \right)  \\
&\cdots  \te_1\big( x_2-x_3+\fr{(n-2)\pi}{n} \bigm| \fr{\tau'}{n} \big) \\
 = \quad
\te_1 \big(n(x_2 - x_3) +(n-2)\pi \mid\tau'\big) & f(x_2)  \te_1\left( x_1-x_3 \bigm| \fr{\tau'}{n} \right)  \te_1\left( x_1-x_3+\fr{\pi}{n} \bigm| \fr{\tau'}{n} \right) \\
& \cdots  \te_1\left( x_1-x_3+\fr{(n-2)\pi}{n} \bigm| \fr{\tau'}{n} \right).
\end{split}
\end{equation*}
{\bf Case 1:\ } $n$ is odd. In this case it is easily seen with the help of (\ref{theta-1}) and (\ref{Transf})
that the function $f(u)= \te_1 \left( u\bigm| \fr{\tau'}{n} \right)$ satisfies
the conditions of Theorem~\ref{master} and with this choice of $f(u)$ the foregoing identity becomes
\begin{equation}\label{help-prod-2}
\begin{split}
&\te_1 \big(n(x_1 -  x_3) +(n-2)\pi\mid\tau'\big)  \te_1 \left( x_1 \bigm| \fr{\tau'}{n} \right)  \te_1\left( x_2-x_3 \bigm| \fr{\tau'}{n} \right) \\
& \qquad \qquad \te_1 \left( x_2-x_3+ \fr{\pi}{n} \bigm| \fr{\tau'}{n} \right)
\cdots  \te_1\big( x_2-x_3+\fr{(n-2)\pi}{n} \bigm| \fr{\tau'}{n} \big) \\
 = & \quad
\te_1 \big(n(x_2 - x_3) +(n-2)\pi\mid\tau'\big)  \te_1 \left( x_2 \bigm| \fr{\tau'}{n} \right)  \te_1\left( x_1-x_3 \bigm| \fr{\tau'}{n} \right)  \\
& \qquad \qquad\te_1\left( x_1-x_3+\fr{\pi}{n} \bigm| \fr{\tau'}{n} \right)
 \cdots  \te_1\left( x_1-x_3+\fr{(n-2)\pi}{n} \bigm| \fr{\tau'}{n} \right).
\end{split}
\end{equation}
Now dividing in (\ref{help-prod-2}) by $x_2-x_3$ and then letting $x_3\to x_2$ implies
\begin{align*}
& \te_1\big(n(x_1-x_2) \mid\tau'\big) \te_1\left( x_1\bigm| \fr{\tau'}{n} \right) \te_1'\left( 0\bigm| \fr{\tau'}{n} \right)
\te_1\left( \fr{\pi}{n}\bigm| \fr{\tau'}{n} \right) \\
& \quad\quad \quad \te_1\left( \fr{2\pi}{n}\bigm| \fr{\tau'}{n} \right)
 \cdots \te_1\left( \fr{(n-2)\pi}{n}\bigm| \fr{\tau'}{n} \right) \\
& = \quad
n \ \te_1'(0\mid\tau') \te_1\left(x_2\bigm| \fr{\tau'}{n}\right)\te_1\left(x_1-x_2\bigm| \fr{\tau'}{n}\right) \te_1\left(x_1-x_2+\fr{\pi}{n}\bigm| \fr{\tau'}{n}\right) \\
& \qquad\qquad\qquad
\te_1\left(x_1-x_2+\fr{2\pi}{n}\bigm| \fr{\tau'}{n}\right) \cdots \te_1\left(x_1-x_2+\fr{(n-2)\pi}{n}\bigm| \fr{\tau'}{n}\right)
\end{align*}
Next, letting $x_1=\fr{\pi}{n}$, this gives
\begin{align*}
& \te_1'\left( 0\bigm| \fr{\tau'}{n} \right) \te_1(nx_2\mid\tau') \te_1^2 \left( \fr{\pi}{n}\bigm| \fr{\tau'}{n} \right) \te_1 \left( \fr{2\pi}{n}\bigm| \fr{\tau'}{n} \right)
\cdots \te_1 \left( \fr{(n-2)\pi}{n}\bigm| \fr{\tau'}{n} \right) \\
& = \quad\quad
- n\ \te_1'(0\mid\tau')\te_1\left(-x_2\bigm| \fr{\tau'}{n}\right) \te_1\left(-x_2+\fr{\pi}{n}\bigm| \fr{\tau'}{n}\right) \\
& \qquad\qquad\qquad
\te_1\left(-x_2+\fr{2\pi}{n}\bigm| \fr{\tau'}{n}\right)
\cdots \te_1\left(-x_2+\fr{(n-1)\pi}{n}\bigm| \fr{\tau'}{n}\right)
\end{align*}
which by substituting $z=-x_2$ and rearranging yields
\begin{align}\label{help-prod-3}
\te_1\left( z\bigm| \fr{\tau'}{n}\right) \te_1\left( z+\fr{\pi}{n}\bigm| \fr{\tau'}{n}\right) \cdots \te_1\left( z+\fr{(n-1)\pi}{n}\bigm| \fr{\tau'}{n}\right)  \\
=
\fr{\te_1'\left(0\bigm|\fr{\tau'}{n}\right)}{n \te_1'(0\mid \tau')} \te_1^2\left(\fr{\pi}{n}\bigm|\fr{\tau'}{n}\right) \te_1\left(\fr{2\pi}{n}\bigm|\fr{\tau'}{n}\right)
\cdots \te_1\left(\fr{(n-2)\pi}{n}\bigm|\fr{\tau'}{n}\right) \te_1(nz\mid\tau'). \nonumber
\end{align}
Thus by virtue of the identities (\ref{Equiv-SineProd}) and (\ref{help-prod-3}) we will be done if we show that
\[
q^{\frac{(n-1)(n+1)}{2}} \frac{(q;q^2)_{\infty}^2}{(q^n;q^{2n})_{\infty}^{2n}}
\frac{\theta_1^n \big( \frac{\pi}{2}\mid \frac{\tau'}{n} \big)}{\theta_1 \big(\frac{\pi}{2}\mid \tau')}
\]
\[
=
\fr{\te_1'\left(0\bigm|\fr{\tau'}{n}\right)}{n \te_1'(0\mid \tau')} \te_1^2\left(\fr{\pi}{n}\bigm|\fr{\tau'}{n}\right) \te_1\left(\fr{2\pi}{n}\bigm|\fr{\tau'}{n}\right)
\cdots \te_1\left(\fr{(n-2)\pi}{n}\bigm|\fr{\tau'}{n}\right),
\]
or, equivalently,
\[
q^{\frac{(n-1)(n+1)}{12}} \frac{(q;q^2)_{\infty}^2}{(q^n;q^{2n})_{\infty}^{2n}}
\]
\begin{equation}\label{help-prod-4}
=
\fr{1}{n}\fr{\te_1'\left(0\bigm|\fr{\tau'}{n}\right)}
{\te_1'(0\mid \tau')}
\fr{\theta_1 \big(\frac{\pi}{2}\mid \tau')}{\theta_1 \big( \frac{\pi}{2}\mid \frac{\tau'}{n} \big)}
\fr{\te_1^2\left(\fr{\pi}{n}\bigm|\fr{\tau'}{n}\right) \prod_{j=2}^{n-2}\te_1\left(\fr{j\pi}{n}\bigm| \fr{\tau'}{n}\right)}
{\theta_1^{n-1} \big( \frac{\pi}{2}\mid \frac{\tau'}{n} \big)}.
\end{equation}
To establish (\ref{help-prod-4}), we use Jacobi's imaginary transformation~(\ref{ImTransf}) and
the infinite product representation~(\ref{theta-product}) and proceed
as follows. For all $j=1,\ldots, n-2$ we have
\[
\begin{split}
\te_1\left(\fr{j\pi}{n}\bigm| \fr{\tau'}{n}\right)
&=
\left(-i \fr{\tau'}{n} \right)^{-\fr{1}{2}}(-i) e^{\fr{i j^2\pi\tau}{n}} \te_1(j\pi\tau \mid n\tau) \\
&=
\left(-i \fr{\tau'}{n} \right)^{-\fr{1}{2}}(-i) e^{\fr{i j^2\pi\tau}{n}}
i q^{\fr{n}{4}} e^{-i j\pi\tau}(q^{2n} e^{-2j i \pi\tau}, e^{2j i \pi\tau}, q^{2n}; q^{2n})_{\infty} \\
&=
\left(-i \fr{\tau'}{n} \right)^{-\fr{1}{2}} q^{\fr{j^2}{n}+\fr{n}{4}-j}(q^{2j},q^{2n-2j},q^{2n};q^{2n})_{\infty},
\end{split}
\]
from which we get
\begin{equation*}
\begin{split}
\te_1^2\left(\fr{\pi}{n}\bigm| \fr{\tau'}{n}\right)
\prod_{j=2}^{n-2}\te_1\left(\fr{j\pi}{n}\bigm| \fr{\tau'}{n}\right)
&=
\left(-i \fr{\tau'}{n} \right)^{-\fr{n-1}{2}} q^N
\prod_{j=1}^{\fr{n-1}{2}}(q^{2j},q^{2n-2j},q^{2n};q^{2n})_{\infty}^2 \\
&=
\left(-i \fr{\tau'}{n} \right)^{-\fr{n-1}{2}} q^N \fr{(q^2;q^2)_{\infty}^2}{(q^{2n};q^{2n})_{\infty}^2}
(q^{2n};q^{2n})_{\infty}^{n-1} \\
&=
\left(-i \fr{\tau'}{n} \right)^{-\fr{n-1}{2}} q^N (q^2;q^2)_{\infty}^2 (q^{2n};q^{2n})_{\infty}^{n-3},
\end{split}
\end{equation*}
where
\[
N= \fr{1}{n}+\fr{n}{4}-1 + \fr{1+2^2+\ldots+ (n-2)^2}{n} + \fr{n(n-2)}{4}-\fr{(n-2)(n-1)}{2}
= \fr{(n-2)(n-1)}{12}.
\]
Similarly,
\[
\begin{split}
\te_1\left(\fr{\pi}{2}\bigm| \fr{\tau'}{n}\right)
&=
\left(-i \fr{\tau'}{n} \right)^{-\fr{1}{2}}(-i) e^{\fr{i n\pi\tau}{4}} \te_1\left(\fr{n\pi\tau}{2} \mid n\tau \right) \\
&=
\left(-i \fr{\tau'}{n} \right)^{-\fr{1}{2}}(-i) e^{\fr{i n\pi\tau}{4}}
i q^{\fr{n}{4}} e^{-\fr{i n\pi\tau}{2}} (q^{2n} e^{-i n \pi\tau}, e^{i n \pi\tau}, q^{2n};q^{2n})_{\infty} \\
&=
\left(-i \fr{\tau'}{n} \right)^{-\fr{1}{2}} (q^n,q^n,q^{2n};q^{2n})_{\infty},
\end{split}
\]
from which we derive
\begin{equation*}
\te_1^{n-1}\left(\fr{\pi}{2}\bigm| \fr{\tau'}{n}\right) = \left(-i \fr{\tau'}{n} \right)^{-\fr{n-1}{2}}
(q^n;q^{2n})_{\infty}^{2n-2} (q^{2n};q^{2n})_{\infty}^{n-1}.
\end{equation*}
From the above we have
\begin{equation}\label{help-prod-5}
\fr{\te_1^2\left(\fr{\pi}{n}\bigm| \fr{\tau'}{n}\right)
\prod_{j=2}^{n-2}\te_1\left(\fr{j\pi}{n}\bigm| \fr{\tau'}{n}\right)}
{\te_1^{n-1}\left(\fr{\pi}{2}\bigm| \fr{\tau'}{n}\right)} =
q^{\fr{(n-2)(n-1)}{12}} \fr{(q^2;q^2)_{\infty}^2}{(q^{2n};q^{2n})_{\infty}^2 (q^n;q^{2n})_{\infty}^{2n-2}}.
\end{equation}
Furthermore, with the help of (\ref{key-derivative}) we find
\begin{align}\label{help-prod-6}
\fr{1}{n}\fr{\te_1'\left(0\bigm|\fr{\tau'}{n}\right)}
{\te_1'(0\mid \tau')}
\fr{\theta_1 \big(\frac{\pi}{2}\mid \tau')}{\theta_1 \big( \frac{\pi}{2}\mid \frac{\tau'}{n} \big)}
&=
\fr{1}{n} \fr{2 (-i n\tau)^{\fr{3}{2}} q^{\fr{n}{4}}(q^{2n};q^{2n})_{\infty}^3
(-i\tau')^{-\fr{1}{2}} (q;q^2)_{\infty}^2 (q^2;q^2)_{\infty}}
{2 (-i \tau)^{\fr{3}{2}} q^{\fr{1}{4}}(q^{2};q^{2})_{\infty}^3
\big(-i\fr{\tau'}{n} \big)^{-\fr{1}{2}} (q^n;q^{2n})_{\infty}^2 (q^{2n};q^{2n})_{\infty}} \nonumber \\
&=
q^{\fr{n-1}{4}} \fr{(q^{2n};q^{2n})_{\infty}^2 (q;q^2)_{\infty}^2}{(q^n;q^{2n})_{\infty}^2 (q^2;q^2)_{\infty}^2}.
\end{align}
Finally, multiply  (\ref{help-prod-5}) and (\ref{help-prod-6}) and simplify to deduce the desired formula
(\ref{help-prod-4}).
\\
{\bf Case 2:\ } If $n$ is even, take $f(u)= \te_1\left(u+\fr{\pi}{2}\bigm| \fr{\tau'}{n}\right)$ and
proceed in exactly the same way to derive the result.
\section{Proof of Theorem~\ref{Gauss-q-Gamma}}
We start proving that
(\ref{q-anlog-Gauss}) is the $q$-analogue for identity~(\ref{Gauss-Gamma-special}). Assuming (\ref{q-anlog-Gauss})
and the basic fact that
$\lim_{q\to 1} \Gamma_q(1/2) = \Gamma (1/2) = \sqrt{\pi}$, it will be
enough to show that
\begin{equation}\label{key-limit}
\lim_{q\to 1}\fr{(q;q^2)_{\infty}^{n-1} (q^2;q^2)_{\infty}}{(q^{\fr{2}{n}};q^{\fr{2}{n}})_{\infty}} =
\fr{2^{\fr{n-1}{2}}}{\sqrt{n}}.
\end{equation}
Note that from (\ref{SineProd}) we have
\[
\prod_{k=1}^{n-1}\sin_{q^n}\pi \left(z+\fr{k}{n} \right) = q^{\frac{(n-1)(n+1)}{12}} \frac{(q;q^2)_{\infty}^2}{(q^n;q^{2n})_{\infty}^{2n}} \fr{\sin_q n\pi z}{\sin_{q^n}\pi z}.
\]
Taking limits as $z\to 0$ on both sides and using (\ref{q-sin-derive}) give
\begin{equation}\label{help-gam-prod-1}
\begin{split}
\prod_{k=1}^{n-1} \sin_{q^n} \fr{k\pi}{n}
&=
 q^{\frac{(n-1)(n+1)}{12}} \frac{(q;q^2)_{\infty}^2}{(q^n;q^{2n})_{\infty}^{2n}}
 \fr{n\pi \sin_q' 0}{\pi \sin_{q^n}'0} \\
&=
q^{\frac{(n-1)(n+1)}{12}} \frac{(q;q^2)_{\infty}^2}{(q^n;q^{2n})_{\infty}^{2n}}  q^{-\fr{n-1}{4}}
\fr{ (q^2;q^2)_{\infty}^2 (q^n;q^{2n})_{\infty}^2}{ (q;q^2)_{\infty}^2 (q^{2n};q^{2n})_{\infty}^2} \\
&=
q^{\fr{(n-1)(n-2)}{12}} \fr{ (q^2;q^2)_{\infty}^2 }{ (q^n;q^{2n})_{\infty}^{2n-2}  (q^{2n};q^{2n})_{\infty}^2} .
\end{split}
\end{equation}
Now let in (\ref{help-gam-prod-1}) $q=q^{\fr{1}{n}}$, next take limits as $q\to 1$, and finally use
the well-known trigonometric formula
\[
\prod_{k=1}^{n-1} \sin \fr{k\pi}{n} = \fr{n}{2^{n-1}}.
\]
to deduce that
\begin{equation*}
\lim_{q\to 1} \fr{(q^{\fr{2}{n}}; q^{\fr{2}{n}})_{\infty}^2}{(q;q^2)_{\infty}^{2n-2} (q^2;q^2)_{\infty}^2}
= \fr{n}{2^{n-1}}
\end{equation*}
which implies (\ref{key-limit}).

\noindent
We now establish the formula (\ref{q-anlog-Gauss}).
By (\ref{sine-cosine-q-gamma}) and (\ref{help-gam-prod-1}) we get
\[
\prod_{k=1}^{n-1} q^{\fr{n}{4}} \Gamma_{q^{2n}}^2\left(\fr{1}{2}\right)
\fr{ q^{n\Big(\fr{k}{n}\big(\fr{k}{n}-1\big) \Big)} }
{ \Gamma_{q^{2n}} \left(\fr{k}{n}\right) \Gamma_{q^{2n}} \left(1-\fr{k}{n}\right) }
=
q^{\fr{(n-1)(n-2)}{12}} \fr{ (q^2;q^2)_{\infty}^2 }{ (q^n;q^{2n})_{\infty}^{2n-2}  (q^{2n};q^{2n})_{\infty}^2},
\]
which after rearranging and simplifying means
\[
 \left( \Gamma_{q^{2n}} \Big(\fr{1}{2}\Big)\right)^{2n-2}
=
 \fr{ (q^2;q^2)_{\infty}^2 }{ (q^n;q^{2n})_{\infty}^{2n-2}  (q^{2n};q^{2n})_{\infty}^2}
\prod_{k=1}^{n-1}\Gamma_{q^{2n}} \left(\fr{k}{n}\right) \Gamma_{q^{2n}} \left(1-\fr{k}{n}\right),
\]
or equivalently,
\[
\left(\prod_{k=1}^{n-1}\Gamma_{q^{2n}}  \left(\fr{k}{n}\right) \right)^2
=
 \left( \Gamma_{q^{2n}} \Big(\fr{1}{2}\Big)\right)^{2n-2}
\fr{ (q^n;q^{2n})_{\infty}^{2n-2}  (q^{2n};q^{2n})_{\infty}^2}{ (q^2;q^2)_{\infty}^2 }.
\]
Now let in the foregoing identity $q = q^{\fr{1}{2n}}$ to obtain
\[
\prod_{k=1}^{n-1}\Gamma_{q}  \Big(\fr{k}{n}\Big)
= \left( \Gamma_{q} \Big(\fr{1}{2}\Big)\right)^{n-1}
\fr{ (q^{\fr{1}{2}};q)_{\infty}^{n-1}  (q;q)_{\infty}}{ (q^{\fr{1}{n}};q^{\fr{1}{n}})_{\infty} }.
\]
This completes the proof.
\section{Proof of Theorem~\ref{q-short-prod}}
By an appeal to Theorem~\ref{Gauss-q-Gamma} and the M\"obius inversion formula, we have
\[
\begin{split}
P_q(n)
&= \prod_{d\mid n}\left(\left(\Gamma_q\left(\fr{1}{2}\right)\right)^{d-1}
\fr{(q^{\fr{1}{2}};q)_{\infty}^{d-1} (q;q)_{\infty}}{(q^{\fr{1}{d}}; q^{\fr{1}{d}})_{\infty}} \right)^{\mu\left(\fr{n}{d}\right)}\\
&= (q;q)_{\infty}^{\sum_{d\mid n}\mu\left(\fr{n}{d}\right)}
\fr{\left(\Gamma_q\left(\fr{1}{2}\right)\right)^{\sum_{d\mid n}d \mu\left(\fr{n}{d}\right)- \sum_{d\mid n}\mu\left(\fr{n}{d}\right)}
(q^{\fr{1}{2}};q)_{\infty}^{\sum_{d\mid n}d\mu\left(\fr{n}{d}\right)}}
{ \prod_{d\mid n}(q^{\fr{1}{d}}; q^{\fr{1}{d}})_{\infty}^{\mu\left(\fr{n}{d}\right)}},
\end{split}
\]
which with the help of the basic facts
\[
\sum_{d\mid n}\mu\left(\fr{n}{d}\right)=0\quad \text{and\quad}  \sum_{d\mid n}d \mu\left(\fr{n}{d}\right) = \varphi(n) \quad (n>1)
\]
gives the desired formula.

\section{Proof of Theorem~\ref{2-powers}}
Let $m>1$ be an integer. The first identity is clear from the definition. As to the second identity, we have
by Theorem~\ref{q-short-prod}, (\ref{q-gamma-half}), and (\ref{q-basics})
\[
\begin{split}
P_q(2^m)
&= \Gamma_{q}^{2^{m-1}}\Big(\fr{1}{2}\Big) \fr{(q^{\fr{1}{2}};q)_{\infty}^{2^{m-1}} (q^{\fr{1}{2^{m-1}}};q^{\fr{1}{2^{m-1}}})_{\infty}}{(q^{\fr{1}{2^m}};q^{\fr{1}{2^m}})_{\infty}} \\
&= (1-q)^{2^{m-2}} \fr{(q;q)_{\infty}^{2^{m-1}} (q^{\fr{1}{2^{m-1}}};q^{\fr{1}{2^{m-1}}})_{\infty}}{(q^{\fr{1}{2^m}};q^{\fr{1}{2^m}})_{\infty}} \\
&= (1-q)^{2^{m-2}} \fr{(q;q)_{\infty}^{2^{m-1}}}{(q^{\fr{1}{2^m}};q^{\fr{1}{2^{m-1}}})_{\infty}}.
\end{split}
\]
Then we will be done if show that
\begin{equation*}
\fr{(q;q)_{\infty}^{2^{m-1}}}{(q^{\fr{1}{2^m}};q^{\fr{1}{2^{m-1}}})_{\infty}} = \psi(q^{\fr{1}{2^{m}}}) \prod_{k=1}^{m-1} \psi^{2^{m-1-k}}(q^{\fr{1}{2^{k}}}).
\end{equation*}
We proceed by induction on $m>1$. If $m=2$, then
\[
\begin{split}
\fr{(q;q)_{\infty}^2}{(q^{\fr{1}{4}};q^{\fr{1}{2}})_{\infty}}
&= \fr{(q;q)_{\infty} (-q^{\fr{1}{2}};q^{\fr{1}{2}})_{\infty} (q^{\fr{1}{2}};q^{\fr{1}{2}})_{\infty}}{(q^{\fr{1}{4}};q^{\fr{1}{2}})_{\infty}} \\
&= \fr{(q;q)_{\infty}}{(q^{\fr{1}{2}};q)_{\infty}} \fr{(q^{\fr{1}{2}};q^{\fr{1}{2}})_{\infty}}{(q^{\fr{1}{4}};q^{\fr{1}{2}})_{\infty}} \\
&= \psi(q^{\fr{1}{4}}) \psi(q^{\fr{1}{2}}),
\end{split}
\]
as required for the basic case. Now suppose the induction hypothesis holds for $m>1$. Then
\[
\begin{split}
\fr{(q;q)_{\infty}^{2^{m}}}{(q^{\fr{1}{2^{m+1}}};q^{\fr{1}{2^{m}}})_{\infty}}
&= (q;q)_{\infty}^{2^{m-1}} \fr{(q;q)_{\infty}^{2^{m-1}}}{ (q^{\fr{1}{2^{m+1}}};q^{\fr{1}{2^{m}}})_{\infty}} \\
&= \fr{(q;q)_{\infty}^{2^{m-1}}}{(q^{\fr{1}{2}};q)_{\infty}^{2^{m-1}}} \fr{(q^{\fr{1}{2}};q^{\fr{1}{2}})_{\infty}^{2^{m-1}}}{\big((q^{\fr{1}{2}})^{\fr{1}{2^{m}}};(q^{\fr{1}{2}})^{\fr{1}{2^{m-1}}}\big)_{\infty}} \\
&= \psi^{2^{m-1}}(q^{\fr{1}{2}}) \psi(q^{\fr{1}{2^{m+1}}}) \prod_{k=1}^{m-1} \psi^{2^{m-1-k}}(q^{\fr{1}{2^{k+1}}}) \\
&= \psi(q^{\fr{1}{2^{m+1}}}) \prod_{k=1}^{m} \psi^{2^{m-k}}(q^{\fr{1}{2^{k}}}).
\end{split}
\]
This completes the proof.

\noindent
{\bf Acknowledgment.} The authors are grateful to the referee for valuable comments and interesting suggestions.
\end{document}